\documentclass[a4paper]{article}

\usepackage{amsmath}
\usepackage{amssymb}
\usepackage[francais]{babel}                                                      
\usepackage{graphicx}

\usepackage[applemac]{inputenc}                                             

\numberwithin{equation}{section}

\newtheorem{thm}{Th\'eor\`eme}[section]                                                              
\newtheorem{lem}[thm]{Lemme}

\newcommand{\ron}[1]{{\cal #1}}

\begin{document}

\begin{center}
{\bf\large Une \'etude asymptotique probabiliste des coefficients d'une s\'erie
enti\`ere}
\end{center}

\begin{center}

Bernard Candelpergher, Michel Miniconi 
\par\smallskip
{\it Laboratoire J.-A. Dieudonn\'e, UMR 7351} 
\par
{\it Universit\'e Nice Sophia Antipolis, 06108 Nice Cedex 02, France} 
\par
{\tt candel@unice.fr, miniconi@unice.fr}
\par\smallskip
(15 juillet 2013)
\end{center}

\medskip

{\it R\'esum\'e. \/} En partant des id\'ees de Rosenbloom~\cite{Rosen} et Hayman~\cite{Hayman}, Luis B\' aez-Duarte donne 
dans~\cite{LBD} une preuve probabiliste de la formule asymptotique de Hardy-Ramanujan pour les partitions d'un entier. 
Le principe g\'en\'eral de la m\'ethode repose sur la convergence en loi d'une famille de variables al\'eatoires vers la loi normale.  Dans notre travail nous d\'emontrons un th\'eor\`eme de type Liapounov (Chung~\cite{Chung}) qui justifie  cette convergence.
L'obtention de  formules asymptotiques simples n\'ecessite une condition dite {\it Gaussienne forte\/} \'enonc\'ee par Luis B\' aez-Duarte, que nous d\'emontrons dans une situation  permettant d'obtenir  une formule asymptotique classique pour les partitions d'un entier en entiers distincts (Erd\" os-Lehner~\cite{Erdos}, Ingham~\cite{Ingham}).\par

\bigskip

{\it Abstract. \/} 
Following the ideas of Rosenbloom~\cite{Rosen} and Hayman~\cite{Hayman}, Luis B\' aez-Duarte gives 
in~\cite{LBD} a probabilistic proof of Hardy-Ramanujan's asymptotic formula for the partitions of an integer. 
The main principle of the method relies on the convergence in law of a family of random variables to a gaussian variable.  In our work we prove a theorem of the Liapounov type (Chung~\cite{Chung}) that justifies  this convergence.
To obtain simple asymptotic formul\ae\ a condition of the so-called {\it strong Gaussian \/} type defined by Luis B\' aez-Duarte is required; we demonstrate this in a situation  that make it possible to obtain a classical asymptotic formula for the partitions of an integer with distinct parts (Erd\"os-Lehner~\cite{Erdos}, Ingham~\cite{Ingham}).

\section{Introduction} 

{\bf Notation.} 

On d\'esigne par $\ron O_{+}(D(0,1))$ l'ensemble des 
fonctions $f$ analytiques de rayon de convergence $1$ telles que $f(t)=\sum_{n\geq 0}a_{n}t^{n}$ avec 
$a_{n}$ r\'eels positifs non tous nuls. En particulier on a $f(t)>0$ pour tout $t\in\ ]0,1[$.
\smallskip

\textbf{D\'efinition.} 

Soit $f\in\ \ron O_{+}(D(0,1))$ : pour tout $t\in\ ]0,1[$ on
d\'efinit la mesure discr\`ete sur $\mathbb{R}$ 
\begin{equation*}
\mu _{t}(f)=\sum_{n\geq 0}\frac{a_{n}t^{n}}{f(t)}\delta _{n}
\end{equation*}

On associe \`{a} cette mesure une variable al\'eatoire $X_{t}$ d\'efinie sur l'espace
probabilis\'e $\Omega =]0,1[,$ \`a valeurs dans $\mathbb{N},$ telle que 
\begin{equation*}
P(X_{t}=n)=\frac{a_{n}t^{n}}{f(t)}
\end{equation*}
le but \'etant de d\'emontrer des r\'esultats sur le comportement asymptotique des $a_n$ en utilisant des m\'ethodes probabilistes ({\it voir\/} Rosenbloom~\cite{Rosen} qui attribue cette id\'ee \`a Khinchin).
\smallskip

{\bf Moments et fonction caract\'eristique.}

Les s\'eries $\sum_{n\geq 0}n^{k}a_{n}t^{n}$ \'etant aussi convergentes
dans le disque $D(0,1),$ on en d\'eduit que $X_{t}$ poss\`ede un moment
d'ordre $k$ pour tout $k\geq 1.$ En particulier pour $k=1$ on a, pour tout $t\in\ ]0,1[$ :
\begin{equation*}
E(X_{t})=\sum_{n\geq 0}n\frac{a_{n}t^{n}}{f(t)}=t\frac{f^{\prime }(t)}{f(t)}
\end{equation*}
On pose $m(t)=E(X_{t})$, $0<t<1$ ; la fonction $m$ est continue sur $]0,1[$.
\par
La fonction caract\'eristique $\varphi _{X_{t}}$ (ou $\varphi _{\mu _{t}(f)}$) de $X_{t}$ est donn\'ee par 
\begin{equation*}
\varphi _{X_{t}}(\theta )=E(e^{i\theta X_{t}})=\sum_{n\geq
0}e^{in\theta }\frac{a_{n}t^{n}}{f(t)}=\frac{f(te^{i\theta })}{f(t)}
\end{equation*}

Les coefficients du d\'eveloppement de Taylor de $f$ sont li\'es \`{a}
la fonction caract\'eristique $\varphi _{X_t}$ par 
\begin{equation} \label{coeffTaylor}
a_{n}=\frac{f(t)}{2\pi t^{n}}\int_{-\pi }^{\pi }\varphi _{X_{t}}(\theta
)e^{-in\theta }d\theta
\end{equation}
pour tout $t\in\ ]0,1[$.

\vfill\eject

{\bf Normalisation.}
\par\smallskip
Soit $\sigma (t)=\sqrt{Var(X_{t})}$ et consid\'erons la variable al\'eatoire
centr\'ee r\'eduite 
\begin{equation*}
Z_{t}=\frac{X_{t}-m(t)}{\sigma (t)}
\end{equation*}

On a 
\begin{equation*}
\varphi _{Z_{t}}(x)=e^{-ix\frac{m(t)}{\sigma (t)}}\varphi _{X_{t}}\bigl(\frac{x}{\sigma (t)}\bigr)
\end{equation*}

donc en posant $\theta =\displaystyle\frac{x}{\sigma (t)}$ dans la formule (\ref{coeffTaylor}) donnant 
$a_{n}$, on obtient pour tout $t\in\ ]0,1 \lbrack $ 
\begin{eqnarray*}
a_{n} =\frac{f(t)}{2\pi \sigma (t)t^{n}}\int_{-\pi \sigma (t)}^{\pi \sigma
(t)}\varphi _{Z_{t}}(x)e^{i\frac{x}{\sigma (t)}(m(t)-n)}dx.
\end{eqnarray*}

Supposons qu'il existe une suite $t_{n}\rightarrow 1$ telle que 
\begin{equation*}
m(t_{n})=n\text{ pour tout }n
\end{equation*}

alors on a 
\begin{equation*}
a_{n}=\frac{f(t_{n})}{2\pi \sigma (t_{n})(t_{n})^{n}}\int_{-\pi \sigma
(t_{n})}^{\pi \sigma (t_{n})}\varphi _{Z_{t_{n}}}(x)dx
\end{equation*}

\bigskip

\textbf{Comportement asymptotique des }$a_{n}$ \textbf{lorsque} $n\rightarrow +\infty $
\par\smallskip
On suppose que la fonction $m$ est continue, strictement croissante et qu'elle tend vers $+\infty$ lorsque $t$ tend vers 1.
Soit  $(t_n)$ une suite dans $]0,1[$ tendant vers 1 et telle que

\begin{equation*}
m(t_{n}) =n\quad \text{pour tout }n
\end{equation*}
et
\begin{equation*}
\sigma (t_{n}) \rightarrow +\infty
\end{equation*}

On a alors
\begin{equation*}
a_{n}=\frac{f(t_{n})}{2\pi \sigma (t_{n})(t_{n})^{n}}\int_{-\pi \sigma
(t_{n})}^{\pi \sigma (t_{n})}\varphi _{Z_{t_{n}}}(x)dx
\end{equation*}

Supposons en outre que l'on ait la convergence en loi de $Z_{t}$ vers une variable al\'eatoire \ $Z$ de loi N(0,1)
quand $t\rightarrow 1$, ce qui veut dire que 
\begin{equation*}
\varphi _{Z_{t}}(x)\rightarrow \varphi _{Z}(x)=e^{-\frac{x^2}{2}} \text{ pour tout } x\in \mathbb{R}, 
\end{equation*}
alors on peut esp\'erer un r\'esultat du type ({\it voir\/} Hayman~\cite{Hayman})
\begin{equation*}
a_{n}\backsim \frac{f(t_{n})}{2\pi \sigma (t_{n})(t_{n})^{n}}\int_{-\infty}^{+\infty }e^{-\frac{x^2}{2}}
dx=\frac{f(t_{n})}{\sqrt{2\pi} \sigma (t_{n})(t_{n})^{n}}
\end{equation*}

Pour que la formule ci-dessus donne un \'equivalent sous une forme  analytique simple, il faudrait pouvoir r\'esoudre explicitement l'\'equation
\begin{equation*}
m(t_{n})=n.
\end{equation*}
Quand ceci  n'est pas  possible, la strat\'egie consiste alors \`a utiliser un \'equivalent de la fonction $t\mapsto
m(t)$ lorsque $t\rightarrow 1$.

Soient $m_1$ et $\sigma_1$ des \'equivalents de $m$ et $\sigma$ respectivement lorsque $t\rightarrow 1$ :
\begin{eqnarray*}
m(t) &\thicksim &m_{1}(t) \\
\sigma (t) &\thicksim &\sigma _{1}(t)
\end{eqnarray*}

\par\smallskip

Soit alors $(\tau_{n})$ une suite dans $]0,1[$ tendant vers $1$ et telle que pour tout $n$
\begin{equation*}
m_{1}(\tau _{n})=n
\end{equation*}

Posons $Z_{t}^{1}=\displaystyle\frac{X_{t}-m_{1}(t)}{\sigma _{1}(t)}$. On a comme
pr\'ec\'edemment

\begin{equation*}
a_{n}=\frac{f(\tau _{n})}{2\pi \sigma _{1}(\tau _{n})\tau _{n}^{n}}%
\int_{-\sigma _{1}(\tau _{n})\pi }^{\sigma _{1}(\tau _{n})\pi
}E(e^{ixZ_{\tau _{n}}^{1}})dx
\end{equation*}

Sous l'hypoth\`ese 
\begin{equation*}
\frac{m(\tau _{n})-m_{1}(\tau _{n})}{\sigma _{1}(\tau _{n})}
\rightarrow 0\text{ quand }\tau _{n}\rightarrow 1
\end{equation*}
et la condition de {\it convergence forte\/} 
\'enonc\'ee par B\' aez-Duarte dans~\cite{LBD} :

\begin{equation*}
\int_{-\sigma (\tau _{n})\pi }^{\sigma (\tau _{n})\pi }\left| \varphi
_{Z_{\tau _{n}}}(x)-e^{-x^{2}/2}\right| dx\rightarrow 0
\end{equation*}
on peut alors obtenir
\begin{equation*}
a_{n}\backsim \frac{f(\tau_{n})}{\sqrt{2\pi} \sigma_1 (\tau_{n})(\tau_{n})^{n}}
\end{equation*}

\medskip

On va appliquer la m\'ethode que l'on vient de d\'ecrire \`a la fonction 
\begin{equation*}
f(z)=\sum q(n)z^{n}
\end{equation*}

o\`u $q(n)$ est le nombre de partitions restreintes de $n$, c'est-\`a-dire le nombre des d\'ecompositions 
$$n=n_1+\cdots+n_p $$ 
en entiers strictement positifs
{\it diff\'erents les uns des autres} afin d'obtenir
la formule asymptotique des par\-ti\-tions restreintes :
\begin{equation*}
q(n)\thicksim\frac{1}{4}\frac{e^{\frac{\pi \sqrt{n}}{\sqrt{3}}}}{3^{1/4}n^{3/4}}
\end{equation*}
({\it voir\/} par exemple Erd\" os~\cite{Erdos} ou Ingham~\cite{Ingham}).

\section{Variable associ\'ee \`{a} un produit infini}
\subsection{\textbf{Mesure associ\'ee \`{a} un produit }}
 
\begin{lem}
\bf Soient $f_{1}$ et $f_{2}$ deux fonctions dans 
$\ron O_{+}(D(0,1))$, alors le produit $f_{1}f_{2}$ est dans
$\ron O_{+}(D(0,1))$ et 
\begin{equation*}
\mu _{t}(f_{1}f_{2})=\mu _{t}(f_{1})\ast \mu _{t}(f_{2})
\end{equation*}

Plus g\'en\'eralement, soient $f_{1}, f_{2},\ldots$
des fonctions dans $\ron O_{+}(D(0,1))$, alors le produit
$f_{1}...f_{n}$ est dans $\ron O_{+}(D(0,1))$ :
\begin{equation*}
\mu _{t}(f_{1}...f_{n})=\mu _{t}(f_{1})\ast\cdots\ast \mu _{t}(f_{n})
\end{equation*}
\end{lem}

\medskip

{\bf D\'emonstration}

Soient $f_{1}(t)=\sum_{n\geq 0}a_{n}t^{n}$ et $f_{2}(t)=\sum_{n\geq
0}b_{n}t^{n},$ on a 
\begin{equation*}
f_{1}(t)f_{2}(t)=\sum_{n\geq 0}a_{n}t^{n}\sum_{n\geq
0}b_{n}t^{n}=\sum_{n\geq 0}\sum_{k+l=n}a_{k}b_{l}t^{n}
\end{equation*}
donc 
\begin{equation*}
\mu _{t}(f_{1}f_{2})=\sum_{n\geq 0}\frac{\sum_{k+l=n}a_{k}b_{l}t^{n}}
             {f_{1}(t)f_{2}(t)}\delta _{n}=\sum_{n\geq 0}\sum_{k+l=n}
             \frac{a_{k}t^{k}b_{l}t^{l}}{f_{1}(t)f_{2}(t)}\delta _{k}\ast \delta _{l}
             =\mu_{t}(f_{1})\ast \mu _{t}(f_{2})
\end{equation*}

Par r\'ecurrence on a $\mu _{t}(f_{1}...f_{n})=\mu _{t}(f_{1})\ast\cdots\ast
\mu _{t}(f_{n}).$

$\square $

\bigskip

\begin{thm}
\bf 
Soit $(f_{n})$ une suite de fonctions dans $\ron O_{+}(D(0,1))$ telle que le
produit infini $\prod_{k\geq 1}^{+\infty }f_{k}$ converge uniform\'ement
sur tout compact de $D(0,1)$. Alors la fonction $f=\prod_{k\geq 1}^{+\infty }f_{k}$
est dans $\ron O_{+}(D(0,1))$ et 
la suite des mesures $\mu _{t}(f_{1}...f_{n})=\mu _{t}(f_{1})\ast\cdots\ast \mu _{t}(f_{n})$ converge en loi vers la mesure $\mu
_{t}(f)$ lorsque $n\rightarrow +\infty $.
\end{thm}

\bigskip

{\bf D\'emonstration.}
Comme le 
produit infini $\prod_{k\geq 1}^{+\infty }f_{k}$ converge uniform\'ement
sur tout compact de $D(0,1)$ on en d\'eduit que la fonction $f=\prod_{k\geq 1}^{+\infty }f_{k}$
est  analytique dans $D(0,1)$.

 En outre on a 
\begin{equation*}
f(z)=\prod_{k\geq 1}^{+\infty }f_{k}(z)=\prod_{k\geq 1}^{+\infty}\sum_{n\geq 0}a_{n,k}z^{n}
=\sum_{n\geq 0}z^{n}\!\!\!\sum_{n_{1}+...+n_{p}=n}a_{n_{1},1}...a_{n_{p},p}
\end{equation*}

donc $f(z)=\sum_{n\geq 0}a_{n}z^{n}$ avec 
\begin{equation*}
a_{n}=\sum_{n_{1}+...+n_{p}=n}a_{n_{1},1}...a_{n_{p},p}
\end{equation*}

ce qui prouve que $f\in\ \ron O_{+}(D(0,1)).$

D'autre part, la fonction caract\'eristique de $\mu _{t}(f_{1}...f_{n})=\mu
_{t}(f_{1})\ast\cdots\ast \mu _{t}(f_{n})$ est \'egale au produit des
fonctions caract\'eristiques de chacune des lois 
\begin{equation*}
\varphi _{\mu _{t}(f_{1})\ast\cdots\ast \mu _{t}(f_{n})}(x)=\frac{f_{1}(te^{ix})}{f_{1}(t)}\cdots\frac{f_{n}(te^{ix})}{f_{n}(t)}
\end{equation*}

Comme le produit $f_{1}(z)...f_{n}(z)$ tend vers $f(z)$ pour tout $z$ dans $D(0,1)$, 
on a pour tout $t\in\ ]0,1[$ 
\begin{equation*}
\lim_{n\rightarrow +\infty }\frac{f_{1}(te^{ix})}{f_{1}(t)}\cdots\frac{f_{n}(te^{ix})}{f_{n}(t)}=\frac{f(te^{ix})}{f(t)}
\end{equation*}

La suite des fonctions caract\'eristiques des mesures $\mu
_{t}(f_{1}...f_{n})$ converge donc simplement vers la fonction
caract\'eristique de la mesure $\mu _{t}(f)$ associ\'ee \`{a} $f.$

$\square $

\bigskip

\subsection{\textbf{La s\'erie des variables al\'eatoires associ\'ees}}
\label{VAssoc}
Soit $(f_{n})$ une suite de fonctions dans $\ron O_{+}(D(0,1))$ telle que le
produit infini $f=\prod_{k\geq 1}^{+\infty }f_{k}$ converge uniform\'ement
sur tout compact de $D(0,1)$.

Par un th\'eor\`eme classique, \`{a} la suite des mesures de
probabilit\'e $(\mu _{t}(f_{n}))$ on peut associer une
probabilit\'e sur l'espace produit $\Omega =]0,1[^{\mathbb{N}}$ et une
suite de variables al\'eatoires  $(X_{n,t})$ \`a valeurs dans $\mathbb{N}$ \textit{ind\'ependantes}
telle que pour tout $n\geq 1$ la variable al\'eatoire $X_{n,t}$ ait pour loi 
$\mu _{t}(f_{n})$.

On peut alors affirmer que pour tout $n\geq 1,$ la mesure $\mu
_{t}(f_{1}...f_{n})$ est la loi de la somme $X_{1,t}+...+X_{n,t}.$

La convergence de la suite de mesures $\mu _{t}(f_{1}...f_{n})$ se traduit
donc par la convergence en loi de la s\'erie $\sum_{n\geq 1}X_{n,t}.$ La
loi de $\sum_{n\geq 1}X_{n,t}$ n'est autre que $\mu _{t}(f).$

\bigskip

D'apr\`es le th\'eor\`eme de Kolmogorov, si la s\'erie 
\begin{equation*}
\sum_{n\geq 1}\sigma ^{2}(X_{n,t})=\sum_{n\geq 1}Var(X_{n,t}-E(X_{n,t}))
\end{equation*}
est convergente, alors  la s\'erie 
$\sum_{n\geq 1}(X_{n,t}-E(X_{n,t}))$ converge presque s\^{u}rement sur 
$\Omega$. Si l'on suppose en outre que la s\'erie  $\sum_{n\geq 1}E(X_{n,t})$ converge alors la s\'erie $\sum_{n\geq 1}X_{n,t}$ converge presque s\^urement.

Les $X_{n,t}$ \'etant positives et ind\'ependantes, on en d\'eduit par le th\'eor\`eme de Beppo-Levi que

\begin{equation*}
Var\bigl(\sum_{n\geq 1}X_{n,t}\bigr)=\sum_{n\geq 1}Var(X_{n,t})
\end{equation*}

Sous ces hypoth\`eses on peut alors donner la d\'efinition suivante :

\bigskip

{
\bf D\'efinition.

Soit $X_t$ la variable al\'eatoire d\'efinie comme la somme  de la s\'erie $\sum_{n\geq 1}X_{n,t}$
dont la loi sera not\'ee $\mu_t(f)$.
On pose
$$
m(t)=E(X_t)=\sum_{n\geq 1}E(X_{n,t})\quad \hbox{et}\quad \sigma (t)=\Bigl(Var(X_{t})\Bigr)^{1/2}=
\Bigl(\sum_{n\geq 1}Var(X_{n,t})\Bigr)^{1/2}
$$
}

\bigskip

\subsection{Application \`{a} l'\'etude des coefficients}

On se donne une suite de fonctions $f_{n}=\sum_{n\geq 0}a_{k,n}z^{k}$  dans
$\ron O_{+}(D(0,1))$ telle que le produit infini $\prod_{n\geq 1}^{+\infty }f_{n}$
converge uniform\'ement sur tout compact de $D(0,1)$. La fonction 
$f=\prod_{n\geq 1}^{+\infty }f_{n}$ est dans $\ron O_{+}(D(0,1))$ et l'on peut \'ecrire
$f(z)=\sum_{n\geq 0}a_{n}z^{n}$ avec 
\begin{equation*}
a_{n}=\sum_{n_{1}+...+n_{p}=n}a_{n_{1},1}...a_{n_{p},p}
\end{equation*}

Si les s\'eries $\sum_{n\geq 1}E(X_{n,t})$\ et 
$\sum_{n\geq 1}\sigma^2(X_{n,t})$\ sont convergentes, consid\'erons la variable al\'eatoire centr\'ee
r\'eduite 
\begin{equation*}
Z_{t}=\frac{\sum_{n\geq 1}X_{n,t}-m(t)}{\sigma (t)}
\end{equation*}

Supposons que $Z_{t}$ converge en loi quand  $t\rightarrow 1$ vers une variable al\'eatoire $Z$ de loi N(0,1) et soit $(t_{n})$ une suite tendant vers 1 telle que 
\begin{equation*}
m(t_{n})=n\text{ pour tout }n.
\end{equation*}

On a alors
\begin{equation*}
\sum_{n_{1}+...+n_{p}=n}a_{n_{1},1}...a_{n_{p},p}=\frac{f(t_{n})}{2\pi
\sigma (t_{n})(t_{n})^{n}}\int_{-\pi \sigma (t_{n})}^{\pi \sigma
(t_{n})}\varphi _{Z_{t_{n}}}(x)dx
\end{equation*}

Si $
\sigma (t_{n})\rightarrow +\infty \text{ quand }t_{n}\rightarrow 1
$
il est plausible que 
\begin{equation*}
\lim_{t_{n}\rightarrow 1}\int_{-\sigma (t_{n})\pi }^{\sigma (t_{n})\pi }\varphi_{t_{n}}(x)dx 
=  \int_{-\infty }^{+\infty }\lim \varphi_{t_{n}}(x)dx=\int_{-\infty }^{+\infty }e^{-x^{2}/2}dx=\sqrt{2\pi }
\end{equation*}

et on en d\'eduirait ainsi \textit{la formule asymptotique des coefficients} 
\begin{equation}    
\label{formuleAsymptGeneral}
\sum_{n_{1}+...+n_{p}=n}a_{n_{1},1}...a_{n_{p},p}\thicksim 
                              \frac{f(t_{n})}{\sqrt{2\pi }\sigma (t_{n})(t_{n})^{n}}.
\end{equation}

\bigskip

{\bf Passage par des \'equivalents.}
\par
Nous d\'etaillons ici la m\'ethode de Luis B\'aez-Duarte~\cite{LBD}.
\smallskip

Afin d'\'etablir la formule (\ref{formuleAsymptGeneral}) ci-dessus avec $\sigma_1$ \`a la place de $\sigma$ 
(et $\tau_n$ \`a la place de $t_n$), il reste \`{a} montrer que 
\begin{equation*}
\lim_{n\rightarrow +\infty }\int_{-\sigma _{1}(\tau _{n})\pi }^{\sigma
_{1}(\tau _{n})\pi }E(e^{ixZ_{\tau _{n}}^{1}})dx=\sqrt{2\pi }.
\end{equation*}

On remarque que l'on peut \'ecrire
\begin{equation*}
Z_{t}^{1}=\frac{X_{t}-m_{1}(t)}{\sigma _{1}(t)}
                          = Z_{t}\frac{\sigma (t)}{\sigma _{1}(t)}+\varepsilon(t)
\end{equation*}
o\`u $\varepsilon(t)=\frac{m(t)-m_{1}(t)}{\sigma _{1}(t)}.$
On a alors
\begin{eqnarray*}
\int_{-\sigma _{1}(\tau _{n})\pi }^{\sigma _{1}(\tau _{n})\pi
}E(e^{ixZ_{\tau _{n}}^{1}})dx &=&\int_{-\sigma _{1}(\tau _{n})\pi }^{\sigma_{1}(\tau _{n})\pi }
                   E(e^{ix\frac{\sigma (\tau _{n})}{\sigma _{1}(\tau _{n})}
                   Z_{\tau _{n}}})e^{ix\varepsilon (\tau _{n})}dx \\
           &=&\frac{\sigma _{1}(\tau _{n})}{\sigma (\tau _{n})}
                   \int_{-\sigma (\tau_{n})\pi }^{\sigma (\tau _{n})\pi }\varphi _{Z_{\tau _{n}}}(x)
                   e^{ix\frac{\sigma _{1}(\tau _{n})}{\sigma (\tau _{n})}\varepsilon (\tau _{n})}dx
\end{eqnarray*}

Pour justifier le remplacement par des \'equivalents on \'enonce deux hypoth\`eses:
\par\smallskip

{\bf Hypoth\`ese 1.} Supposons que 
\begin{equation*}
\frac{m(\tau _{n})-m_{1}(\tau _{n})}{\sigma _{1}(\tau _{n})}=\varepsilon
(\tau _{n})\rightarrow 0\text{ quand }\tau _{n}\rightarrow 1
\end{equation*}

{\bf Hypoth\`ese 2.} Supposons  que 
\begin{equation*}
\int_{-\sigma (\tau _{n})\pi }^{\sigma (\tau _{n})\pi }\left| \varphi
_{Z_{\tau _{n}}}(x)-e^{-x^{2}/2}\right| dx\rightarrow 0
\end{equation*}

Sous ces deux hypoth\`eses il est facile de montrer que la suite $\displaystyle\int_{-\sigma (\tau _{n})\pi }^{\sigma (\tau _{n})\pi
}\varphi _{Z_{\tau _{n}}}(x)e^{ix\frac{\sigma _{1}(\tau _{n})}{\sigma (\tau
_{n})}\varepsilon (\tau _{n})}dx$ converge vers $\sqrt{2\pi }$. 

En effet,  on a

\begin{eqnarray*}
&&\left| \int_{-\sigma (\tau _{n})\pi }^{\sigma (\tau _{n})\pi }\varphi_{Z_{\tau _{n}}}(x)
           e^{ix\frac{\sigma _{1}(\tau _{n})}{\sigma (\tau _{n})}
          \varepsilon (\tau _{n})}dx
          -\int_{-\sigma (\tau _{n})\pi }^{\sigma (\tau_{n})\pi }
          e^{ix\frac{\sigma _{1}(\tau _{n})}{\sigma (\tau _{n})}\varepsilon
          (\tau _{n})}e^{-x^{2}/2}dx\right| \\
          &\leq &\int_{-\sigma (\tau _{n})\pi }^{\sigma (\tau _{n})\pi }\left| 
          \varphi_{Z_{\tau _{n}}}(x)-e^{-x^{2}/2}\right| dx
\end{eqnarray*}
il suffit donc de montrer que 

\begin{equation*}
\lim_{\tau_n\rightarrow 1}\int_{-\sigma (\tau _{n})\pi }^{\sigma (\tau _{n})\pi }
        e^{ix\frac{\sigma_{1}(\tau _{n})}{\sigma (\tau _{n})}\varepsilon (\tau_{n})}e^{-x^{2}/2}dx = \sqrt{2\pi }.
\end{equation*}

Ceci r\'esulte du th\'eor\`eme de la convergence domin\'ee, car on a 
$$
\lim_{\tau_n\rightarrow 1}e^{ix\frac{\sigma
_{1}(\tau _{n})}{\sigma (\tau _{n})}\varepsilon (\tau
_{n})}e^{-x^{2}/2} = e^{-x^{2}/2}
$$ 
et 
$$
\left| e^{ix\frac{\sigma
_{1}(\tau _{n})}{\sigma (\tau _{n})}\varepsilon (\tau
_{n})}e^{-x^{2}/2}\right| \leq e^{-x^{2}/2}
$$

R\'esumons ce qui pr\'ec\`ede dans le th\'eor\`eme suivant:

\begin{thm}
\label{thmEquiv}
\bf

(Th\'eor\`eme des \'equivalents)
Soit $(f_{n}=\sum_{n\geq 0}a_{k,n}z^{k})$ une suite de fonctions dans 
$\ron O_{+}(D(0,1))$ telle que le produit infini $\prod_{n\geq 1}^{+\infty }f_{n}$
converge uniform\'ement sur tout compact de $D(0,1)$. La fonction 
$f=\prod_{n\geq 1}^{+\infty }f_{n}$ est donc analytique dans $D(0,1)$ et 
$f(z)=\sum_{n\geq 0}a_{n}z^{n}$ avec 
\begin{equation*}
a_{n}=\sum_{n_{1}+...+n_{p}=n}a_{n_{1},1}...a_{n_{p},p}
\end{equation*}

Si les s\'eries $\sum_{n\geq 1}E(X_{n,t})$\ et $\sum_{n\geq 1}\sigma^2(X_{n,t})$
sont convergentes, consid\'erons la variable al\'eatoire centr\'ee r\'eduite 
\begin{equation*}
Z_{t}=\frac{\sum_{n\geq 1}X_{n,t}-m(t)}{\sigma (t)}
\end{equation*}

Supposons que $Z_{t}$ converge en loi quand $t\rightarrow 1$ vers une variable al\'eatoire $Z$ de loi N(0,1) avec 
\textit{la condition de convergence forte} ({\it voir\/}~\cite{LBD}) :
\begin{equation} \label{condConvForte}
\lim_{t\rightarrow 1}\int_{-\sigma (t)\pi }^{\sigma (t)\pi }\left| \varphi
_{Z_{t}}(x)-e^{-x^{2}/2}\right| dx=0
\end{equation}

Soient $m_1$ et $\sigma_1$ des \'equivalents de $m$ et $\sigma$ respectivement lorsque $t\rightarrow 1$ :
\begin{eqnarray*}
m(t) &\thicksim &m_{1}(t) \\
\sigma (t) &\thicksim &\sigma _{1}(t)
\end{eqnarray*}

Soit une suite $(\tau _{n})$ dans $]0,1[$ convergeant vers 1 et telle que pour tout $n$ on ait : 
\begin{equation*}
m_{1}(\tau _{n})=n
\end{equation*}

avec 
\begin{equation*}
\frac{m(\tau _{n})-m_{1}(\tau _{n})}{\sigma _{1}(\tau _{n})}\rightarrow 0\text{ quand }\tau _{n}\rightarrow 1.
\end{equation*}

Alors

\begin{equation}
\label{formuleAsymptGeneralTau}
a_{n}\thicksim \frac{f(\tau _{n})}{\sqrt{2\pi }\sigma _{1}(\tau _{n})\tau_{n}^{n}}
\end{equation}

\end{thm}

\section{Un th\'eor\`eme de convergence}

Nous \'enon\c cons et d\'emontrons un th\'eor\`eme du type Liapounov ({\it voir\/} Chung~\cite{Chung} p. 205 sq.) de convergence vers une loi normale concernant une famille continue de suites infinies de variables al\'eatoires.

\begin{thm}

\label{thmConv}
\bf (Th\'eor\`eme de convergence)
Soit une suite de variables al\'eatoires positives $(X_{n,t})_n$ dans $L^3(\Omega)$
telle que pour tout $t\in\ ]0,1[$ :
\begin{description}
\item{a)} les $X_{n,t}$ sont ind\'ependantes

\medskip       

\item{b)} les s\'eries $m(t) = \sum_{n\geq 1}E(X_{n,t})$, 
$\sigma^2(t) = \sum_{n\geq 1}Var(X_{n,t})$ et 
$\Gamma_{3}(t)=\sum_{n\geq 1}E(\left| X_{n,t}-E(X_{n,t})\right| ^{3})$ sont
convergentes

\medskip

\item{c)} la fonction $t\mapsto\displaystyle \frac{\Gamma_{3}(t)}{(\sigma (t))^{3}}$ tend vers $0$ 
quand $t\rightarrow 1$

\medskip
\item{d)} ${\boldsymbol\lim_{t\rightarrow 1}}{\boldsymbol\sup}_{n\geq 1}\displaystyle\frac{Var(X_{n,t})}{\sigma ^{2}(t)}=0$

\end{description}
\medskip
Alors la s\'erie $X_t=\sum_{n\geq 1}X_{n,t}$ est convergente presque sûrement et
la fonction caract\'eristique $\varphi _{Z_{t}}$ de la variable al\'eatoire 
\begin{equation*}
Z_{t}=\frac{X_t-m(t)}{\sigma (t)}
\end{equation*}

est telle que 
\begin{equation*}
\varphi _{Z_{t}}(x)\rightarrow e^{-x^{2}/2}\quad\hbox{quand}\quad 
t\rightarrow 1.
\end{equation*}

\end{thm}

\bigskip

{\bf D\'emonstration.}

\smallskip
Posons $Y_{n,t}=X_{n,t}-E(X_{n,t})$ on a $\displaystyle Z_{t}=\frac{\sum_{n\geq 1}Y_{n,t}}{\sigma (t)}.$

Par l'ind\'ependance des $Y_{n,t}$ on voit que la variable al\'eatoire $Z_{t}$ a pour
fonction caract\'eristique 
\begin{equation*}
\varphi _{Z_{t}}(\theta )=E(e^{i\theta \sum \frac{Y_{n,t}}{\sigma (t)}})
          = \prod_{n\geq 1}E(e^{i\theta \frac{Y_{n,t}}{\sigma (t)}})
\end{equation*}

On a $E(Y_{n,t})=0$ et $E(Y_{n,t}^{2})=\sigma _{n,t}^{2}=Var(X_{n,t}).$

\begin{lem}
\bf Sous les hypoth\`eses du th\'eor\`eme (\ref{thmConv}) ci-dessus, on a
\begin{equation*}
E\bigl(e^{i\theta \frac{Y_{n,t}}{\sigma (t)}}\bigr)=1-\frac{\theta ^{2}}{2}
               \bigl(\frac{\sigma _{n,t}}{\sigma (t)}\bigr)^{2}+L_{n}(\theta ,t)
\end{equation*}

avec 
\begin{equation*}
\left| L_{n}(\theta ,t)\right| \leq \frac{\left| \theta \right| ^{3}}{6(\sigma (t))^{3}}E(\left| Y_{n,t}\right| ^{3}).
\end{equation*}

\end{lem}
\smallskip

Ce lemme r\'esulte de la formule de Taylor 
\begin{equation*}
e^{ix}=1+ix-\frac{x^{2}}{2}-i\int_{0}^{1}\frac{(1-u)^{2}}{2}x^{3}e^{iux}du
\end{equation*}

qui nous donne\begin{equation*}
e^{i\theta \frac{Y_{n,t}}{\sigma (t)}}=1+i\theta \frac{Y_{n,t}}{\sigma (t)}-
\frac{\theta ^{2}(\frac{Y_{n,t}}{\sigma (t)})^{2}}{2}-i\int_{0}^{1}\frac{(1-u)^{2}}{2}\theta ^{3}(\frac{Y_{n,t}}{\sigma (t)})^{3}
               e^{iu\theta \frac{Y_{n,t}}{\sigma (t)}}du
\end{equation*}

Comme $E(Y_{n,t})=0$ on en d\'eduit que 
\begin{eqnarray*}
E\Bigl(e^{i\theta \frac{Y_{n,t}}{\sigma (t)}}\Bigr) &=&1-\frac{\theta ^{2}
                         (\frac{\sigma _{n,t}}{\sigma (t)})^{2}}{2}-i\int_{0}^{1}\frac{(1-u)^{2}}{2}\theta
^{3}E\Bigl(\bigl(\frac{Y_{n,t}}{\sigma (t)}
                        \bigr)^{3}e^{iu\theta \frac{Y_{n,t}}{\sigma (t)}}\Bigr)du \\
&=&1-\frac{\theta ^{2}}{2}(\frac{\sigma _{n,t}}{\sigma (t)})^{2}+L_{n}(\theta ,t)
\end{eqnarray*}

o\`u 
\begin{equation*}
L_{n}(\theta ,t)=-i\int_{0}^{1}\frac{(1-u)^{2}}{2}\theta ^{3}
              E\Bigl(\bigl(\frac{Y_{n,t}}{\sigma (t)}\bigr)^{3}e^{iu\theta \frac{Y_{n,t}}{\sigma (t)}}\Bigr)du
\end{equation*}

On a ainsi la majoration
\begin{equation*}
\left| L_{n}(\theta ,t)\right| \leq |\theta ^{3}|E\Bigl(\bigl(\frac{|Y_{n,t}|}{\sigma (t)}
\bigr)^{3}\Bigr)\int_{0}^{1}\frac{(1-u)^{2}}{2}\, du\leq \frac{\left| \theta \right| ^{3}}{6(\sigma (t))^{3}}E(\left| Y_{n,t}\right| ^{3})
\end{equation*}

Ce qui termine la d\'emonstration du lemme.

$\square $

\bigskip

Par l'ind\'ependance des $Y_{n,t}$ la fonction
caract\'eristique $\varphi _{Z_{t}}$ de la variable al\'eatoire $Z_{t}$ peut s'\'ecrire  
\begin{equation}
\label{fonctCaract}
\varphi _{Z_{t}}(\theta )=\prod_{n\geq 1}\bigl(1-\frac{\theta ^{2}}{2}
                  \bigl(\frac{\sigma _{n,t}}{\sigma (t)}\bigr)^{2}+L_{n}(\theta ,t)\bigr)
\end{equation}

Pour  montrer que
$\varphi _{Z_{t}}(\theta )\rightarrow e^{-\frac{\theta ^{2}}{2}}$ quand $t\rightarrow 1$ 
nous allons utiliser le lemme suivant (dont nous donnons la d\'emonstration dans l'Appendice ({\it voir\/} section 5)) :

\begin{lem}  
\label{lemmeConv}
\bf Soit $(u_{n,t})_{n\geq 1}$ une famille de suites complexes
index\'ees par $t\in\ ]0,1[$ telle que
\begin{description}
\item{($\imath$)} $\boldsymbol\sup_{n\geq 1}|u_{n,t}|\rightarrow 0$ quand $t\rightarrow 1.$
\smallskip
\item{($\imath\imath$)} il existe $M>0$ et $0<\alpha <1$ tel que $\sum_{n\geq 1}|u_{n,t}|\leq M$
pour tout $t\in\ ]\alpha ,1[.$
\smallskip
\item{($\imath\imath\imath$)} il existe $S\in\ \mathbb{C}$ tel que $\sum_{n\geq 1}u_{n,t}\rightarrow S 
$ quand $t\rightarrow 1.$

\end{description}
\smallskip
Alors 
$$
\lim_{t\rightarrow 1}\prod_{n\geq 1}(1+u_{n,t}) = e^{S}   
$$
\end{lem}

\medskip

On va appliquer ce lemme \`a la fonction caract\'eristique (\ref{fonctCaract}) en posant 
\begin{equation*}
u_{n,t}(\theta )=-\frac{\theta ^{2}}{2}(\frac{\sigma _{n,t}}{\sigma (t)}
)^{2}+L_{n}(\theta ,t)
\end{equation*}

V\'erifions les trois conditions du lemme:

\smallskip
($\imath$) On a 
\begin{eqnarray*}
\sup_{n\geq 1}|u_{n,t}(\theta )| &\leq &\frac{\theta ^{2}}{2}\sup_{n\geq 1}
           (\frac{\sigma _{n,t}}{\sigma (t)})^{2}+\sup_{n\geq 1}L_{n}(\theta ,t) \\
&\leq &\frac{\theta ^{2}}{2}\sup_{n\geq 1}(\frac{\sigma _{n,t}}{\sigma (t)})^{2}+
           \frac{\left| \theta \right| ^{3}}{6(\sigma (t))^{3}}\sum_{n\geq 1}E(\left| Y_{n,t}\right| ^{3})
\end{eqnarray*}

D'apr\`es les hypoth\`eses  c) et d) du th\'eor\`eme cette derni\`ere quantit\'e tend vers $0$ quand $t\rightarrow 1.$

\smallskip
($\imath\imath$) On a 
\begin{eqnarray*}
\sum_{n\geq 1}|u_{n,t}(\theta )| &\leq &\frac{\theta ^{2}}{2}\sum_{n\geq 1}
                 (\frac{\sigma _{n,t}}{\sigma (t)})^{2}+\sum_{n\geq 1}L_{n}(\theta ,t) \\
&\leq &\frac{\theta ^{2}}{2}+\frac{\left| \theta \right| ^{3}}{6(\sigma(t))^{3}}\sum_{n\geq 1}E(\left| Y_{n,t}\right| ^{3})
\end{eqnarray*}

Or d'apr\`es c) la quantit\'e $\frac{1}{(\sigma (t))^{3}}\sum_{n\geq
1}E(\left| Y_{n,t}\right| ^{3})$ est born\'ee au voisinage de 1.

\smallskip
($\imath\imath\imath$) On a 
\begin{equation*}
\sum_{n\geq 1}u_{n,t}(\theta )=-\frac{\theta ^{2}}{2}\sum_{n\geq 1}
                  (\frac{\sigma _{n,t}}{\sigma (t)})^{2}+\sum_{n\geq 1}L_{n}(\theta ,t)=\frac{\theta
^{2}}{2}+\sum_{n\geq 1}L_{n}(\theta ,t)
\end{equation*}
et $\sum_{n\geq 1}L_{n}(\theta ,t)\rightarrow 0$ quand $t\rightarrow 1$ par
c).

\smallskip
On a donc $\sum_{n\geq 1}u_{n,t}(\theta )\rightarrow -\frac{\theta ^{2}}{2}$
quand $t\rightarrow 1$ et par le lemme (\ref{lemmeConv})
\begin{equation*}
\varphi _{Z_{t}}(\theta )=\prod_{n\geq 1}\bigl(1-\frac{\theta ^{2}}{2}
              \bigl(\frac{\sigma _{n,t}}{\sigma (t)}\bigr)^{2}+L_{n}(\theta ,t)\bigr)\rightarrow e^{-\frac{\theta ^{2}}{2}}
\end{equation*}

$\square $

\vfill\eject
 
\section{Application aux partitions restreintes}

\label{partRestr}

Consid\'erons la fonction d\'efinie par le produit infini 
$$f(z)=\prod_{n\geq 1}^{+\infty }(1+z^{n})$$ 
Cette fonction est analytique dans \vrule height 0pt depth 8pt width 0pt
$D(0,1)$ car la s\'erie $\sum_{n=1}^{+\infty }z^{n}$ converge
uniform\'ement sur tout compact de $D(0,1)$. On a 
\begin{equation*}
f(z)=\sum q(n)z^{n}
\end{equation*}

o\`u $q(n)$ est le nombre de partitions restreintes de $n,$ c'est-\`a-dire le nombre des d\'ecompositions 
$n=n_1+\cdots+n_p $ en entiers strictement positifs
{\it diff\'erents les uns des autres.}

\smallskip

Le but de ce qui suit est d'appliquer la m\'ethode d\'ecrite au d\'ebut de cet article pour obtenir la formule asymptotique des partitions restreintes :
\begin{equation*}
q(n)\thicksim\frac{1}{4}\frac{e^{\frac{\pi \sqrt{n}}{\sqrt{3}}}}{3^{1/4}n^{3/4}}
\end{equation*}

Soit la mesure de probabilit\'e associ\'ee \`{a} $f_{n}(t)=1+t^{n}$
\begin{equation*}
\mu _{t}(f_{n})=\frac{1}{1+t^{n}}\delta _{0}+\frac{t^{n}}{1+t^{n}}\delta _{n}
\end{equation*}
o\`u $\delta _0$ et $\delta _n$ repr\'esentent les mesures de Dirac en $0$ et $n$ respectivement.

\smallskip
On associe \`a ces mesures une suite de variables al\'eatoires ind\'ependantes $(X_{n,t})$ ({\it voir} section (\ref{VAssoc})).
La variable $X_{n,t}$ prend les valeurs $0$ et $n$ et on a 
\begin{eqnarray*}
E(X_{n,t}) &=&\frac{nt^{n}}{1+t^{n}} \\
Var(X_{n,t}) &=&\frac{n^{2}t^{n}}{\left( 1+t^{n}\right) ^{2}}
\end{eqnarray*}

\medskip

\textbf{Dans ce qui suit on posera} 
$$
t=e^{-r }
$$
\textbf{o\`u} $r >0,$ \textbf{de sorte que l'on a} 
$$
t\rightarrow 1\Leftrightarrow r \rightarrow 0
$$

\subsection{V\'erification des hypoth\`eses du th\'eor\`eme de convergence}

Les s\'eries 
\begin{equation*}
m(t)=  \sum_{n\geq 1}E(X_{n,t}) = \sum_{n=1}^{+\infty }\frac{nt^{n}}{1+t^{n}}
\quad\hbox{et}\quad
\sigma^{2}(t) = \sum_{n\geq 1}\sigma^2(X_{n,t}) = \sum_{k\geq 1}\frac{n^{2}t^{n}}{\left( 1+t^{n}\right) ^{2}}
\end{equation*}
sont clairement convergentes.

Examinons la s\'erie $\sum_{n\geq 1}E(\left| X_{n,t}-E(X_{n,t})\right| ^{3})$ : on a 
\begin{eqnarray*}
E(\left| X_{n,t}-E(X_{n,t})\right| ^{3}) &=&(\frac{nt^{n}}{1+t^{n}})^{3}
\frac{1}{1+t^{n}}+(n-\frac{nt^{n}}{1+t^{n}})^{3}\frac{t^{n}}{1+t^{n}} \\
&=&n^{3}\frac{t^{3n}+t^{n}}{\left( 1+t^{n}\right) ^{4}}
\end{eqnarray*}

donc la s\'erie $\sum_{n\geq 1}E(\left| X_{n,t}-E(X_{n,t})\right| ^{3})$ est
convergente.

\smallskip
Ainsi les hypoth\`eses a) et b) du th\'eor\`eme de convergence (\ref{thmConv}) sont bien v\'erifi\'ees.

\subsubsection{Comportement asymptotique de $m$ et $\sigma^2$}

Pour d\'eterminer le comportement asymptotique quand $t\rightarrow 1$ des fonctions $m(t)$ et $\sigma^2(t)$ on va utiliser la formule d'Euler-McLaurin rappel\'ee ci-dessous :

si $f\in\ C^{1}[0,+\infty \lbrack $ on a pour tout entier $n\geq 1$
\begin{equation*}
\sum_{k=1}^{n}f(k)=\int_{1}^{n}f(x)dx+\frac{1}{2}(f(1)+f(n))+
\int_{1}^{n}b_{1}(x)f^{\prime }(x)dx
\end{equation*}
o\`u $b_{1}(x)=x-[x]-\frac{1}{2}$. 
Si en outre $\sum_{k=1}^{+\infty }f(k)$ et $\int_{1}^{+\infty }f(x)dx$ sont convergentes alors 
\begin{equation*}
\sum_{k=1}^{+\infty }f(k)=\int_{0}^{+\infty }f(x)dx+\int_{1}^{+\infty }b_{1}(x)f^{\prime }(x)dx+C
\end{equation*}
o\`u  $C=\frac12f(1)-\int_{0}^{1}f(x)dx.$
\bigskip

Les fonctions $f$ auxquelles on va appliquer cette formule seront du type 
\begin{equation*}
f(x)=\frac{x^{p}e^{-arx}}{(1+e^{-rx})^{q}}
\end{equation*}
o\`u $a,p,q$ sont des entiers sup\'erieurs ou \'egaux \`a 1.
Comme $ f(x)=\frac{1}{r^p}g(rx)$ o\`u  $g(u)=\frac{u^{p}e^{-au}}{(1+e^{-u})^{q}}$, on a
\begin{equation*}
\int_{0}^{+\infty }f(x)dx=\frac{1}{r^{p+1}}\int_{0}^{+\infty }g(u)du
\end{equation*}
et
\begin{equation*}
\Bigl|\int_{1}^{+\infty }b_{1}(x)f^{\prime }(x)dx\Bigr| =
\frac{1}{r^{p}}\Bigl|\int_{1}^{+\infty }b_{1}(\frac ur)g^{\prime }(u)du\Bigr|\leq \frac{1}{2r^{p}} \int_{1}^{+\infty }\vert g^{\prime }(u)\vert du
\end{equation*}
car la fonction $g^{\prime }$ est int\'egrable.

\medskip

La formule d'Euler-MacLaurin nous donne pour  $r\rightarrow 0+$
\begin{equation*}
\sum_{k=1}^{+\infty}\frac{k^{p}e^{-akr}}{(1+e^{-kr})^{q}}=\frac{1}{r^{p+1}}\int_{0}^{+\infty }
\frac{u^{p}e^{-au}}{(1+e^{-u})^{q}}+O(\frac{1}{r^p}).
\end{equation*}

Pour $a=p=q=1$ on obtient  
\begin{eqnarray*}
m(e^{-r })=\sum_{k=1}^{+\infty }\frac{ke^{-r k}}{1+e^{-r k}} =\frac{1}{r^2}\int_{0}^{+\infty
}\frac{ue^{-u}}{1+e^u}dx+O(\frac 1r)
\end{eqnarray*}

Notons que

\begin{eqnarray*}
\frac{1}{r^2}\int_{0}^{+\infty
}\frac{ue^{-u}}{1+e^u}dx
=\frac{1}{r^2}\sum_{n=0}^{+\infty }(-1)^{n}\int_{0}^{+\infty }ue^{-u(n+1)}dx
=\frac{1}{r ^{2}}\sum_{n=0}^{+\infty }(-1)^{n}\frac{1}{\left(
n+1\right) ^{2}}=\frac{1}{r ^{2}}\frac{\pi ^{2}}{12}
\end{eqnarray*}

\medskip

On a ainsi
\begin{equation*}
m(e^{-r })=m_{1}(e^{-r })+O(\frac 1r)\quad\hbox{avec}\quad  m_{1}(e^{-r })=\frac{\pi ^{2}}{12}\frac{1}{r ^{2}}
\end{equation*}

\medskip

Et de la m\^eme mani\`ere, on a

\begin{equation}       \label{sigma2}
\sigma ^{2}(e^{-r })\thicksim \int_{1}^{+\infty }
                \frac{x^{2}e^{-r x}}{(1+e^{-r x})^{2}}dx
                \thicksim _{r \rightarrow 0}\frac{\pi ^{2}}{6}\frac{1}{r ^{3}}=\sigma _{1}^{2}(e^{-r })
\end{equation}

car

\begin{equation*}
\int_{0}^{+\infty }\frac{x^{2}e^{-r x}}{(1+e^{-r x})^{2}}dx
                 =\sum_{n=1}^{+\infty }(-1)^{n-1}n\int_{0}^{+\infty }x^{2}e^{-r xn}dx
                 =\frac{2}{r ^{3}}\sum_{n=1}^{+\infty }(-1)^{n-1}\frac{1}{n^{2}}
\end{equation*}

\bigskip

\subsubsection{Les conditions c) et d)}
Calculons $\Gamma _{3}(t)=\sum_{n\geq 1}E(\left| X_{n,t}-E(X_{n,t})\right|^{3})$ :
\begin{equation*}
\sum_{n\geq 1}E(\left| X_{n,t}-E(X_{n,t})\right| ^{3})=\sum_{n\geq 1}n^{3}
\frac{e^{-3nr }+e^{-nr }}{\left( 1+e^{-nr }\right) ^{4}}\thicksim
\int_{0}^{+\infty }\frac{x^{3}(e^{-3r x}+e^{-r x})}{(1+e^{-r x})^{4}}dx=\frac{C}{r ^{4}}
\end{equation*}
et donc 
\begin{equation*}
\frac{\Gamma _{3}(t)}{\sigma (t)^{3}}\thicksim \frac{\frac{C}{r ^{4}}}{(\frac{\pi ^{2}}{6}\frac{1}{r ^{3}})^{3/2}}=C_{3}r ^{1/2}
\end{equation*}

On a donc bien $\frac{\Gamma _{3}(t)}{\sigma (t)^{3}}\rightarrow 0$ quand $t\rightarrow 1.$

Il reste \`{a} voir que 
$\lim_{t\rightarrow 1}\sup_{n\geq 1}\displaystyle\frac{Var(X_{n,t})}{\sigma ^{2}(t)}=0$. On a

\begin{equation*}
\frac{Var(X_{n,t})}{\sigma ^{2}(t)}=\frac{1}{\sigma ^{2}(t)}n^{2}
\frac{t^{n}}{\left( 1+t^{n}\right) ^{2}}\leq \frac{1}{\sigma ^{2}(t)}n^{2}t^{n}
\end{equation*}

Or on a $n^{2}e^{-nr }\leq \frac{4}{r ^{2}}e^{-2}$ pour tout $n$ et 
$\sigma ^{2}(e^{-r })\thicksim \frac{\pi ^{2}}{6}\frac{1}{r ^{3}}$ d'apr\`es (\ref{sigma2}) donc
$$
\lim_{t\rightarrow 1}\sup_{n\geq 1}\frac{Var(X_{n,t})}{\sigma ^{2}(t)}=0
$$

\bigskip

Ainsi les hypoth\`eses du th\'eor\`eme de convergence (\ref{thmConv}) sont  bien v\'erifi\'ees. Par cons\'equent
la fonction caract\'eristique $\varphi _{Z_{t}}$ de la variable al\'eatoire 
\begin{equation*}
Z_{t}=\frac{\sum_{n\geq 1}X_{n,t}-m(t)}{\sigma (t)}
\end{equation*}
converge vers $e^{-x^{2}/2}$ quand $t\rightarrow 1$.

\smallskip

\bigskip

\subsection{V\'erification de la condition de convergence forte}
Pour obtenir une formule asymptotique du nombre de partitions restreintes $q(n)$ d\'efini au d\'ebut de ce paragraphe (\ref{partRestr}), on doit v\'erifier
les hypoth\`eses du th\'eor\`eme des \'equivalents, en particulier
la condition de convergence forte :

\begin{equation*}
\lim_{t\rightarrow 1}\int_{-\pi \sigma (t)}^{\pi \sigma (t)}\left| \varphi _{Z_t}(\theta)-e^{-\theta ^{2}/2}\right| d\theta = 0
\end{equation*}
Pour cela on va d\'ecomposer l'int\'egrale pr\'ec\'edente en 
$$\int_{\left| \theta \right| \leq \frac{C}{r ^{1/2}}}\left| \varphi
_{Z_t}(\theta )-e^{-\theta ^{2}/2}\right| d\theta +
\int_{\frac{C}{r ^{1/2}}\leq \left| \theta \right| \leq \pi \sigma
(t)}\left| \varphi _{Z_t}(\theta )-e^{-\theta ^{2}/2}\right| d\theta  
$$
et majorer $\left| \varphi _{Z_t}(\theta )\right|$ sur chacun des domaines d'int\'egration.

\begin{lem} \label{lemmeA}
\bf  Si $\left| \theta \right| \leq \displaystyle \frac{1}{4\frac{\Gamma
_{3}(t)}{\sigma ^{3}(t)}}\thicksim \frac{1}{4C_{3}}r ^{-1/2}$ 
alors $\left| \varphi _{Z_t}(\theta )\right| \leq e^{-\theta ^{2}/3}$.
\end{lem}

{\bf D\'emonstration.}

Posons $Y_{n,t}=X_{n,t}-E(X_{n,t})$ on a $Z_{t}=\frac{\sum_{n\geq 1}Y_{n,t}}{\sigma (t)}.$ On a 
\begin{equation*}
\varphi _{Z_t}(\theta )=\prod_{n\geq 1}\varphi _{Y_{n,t}}\bigl(\frac{\theta }{\sigma (t)}\bigr)
\end{equation*}

Pour majorer $\left| \varphi _{Z_t}(\theta )\right|  $ on va utiliser  le lemme suivant dont la d\'emonstration est report\'ee \`a la section (\ref{preuveLemCramer}) : 

\begin{lem} \label{lemmeCramer} 
\bf (Lemme de Cram\'er) Soit $Z$\ une variable al\'eatoire centr\'ee telle que $E(\left| Z\right| ^{3})<+\infty $ et $\varphi_{Z}$ sa fonction caract\'{e}ristique. 
On a 
\begin{equation*}
\left| \varphi _{Z}(\xi)\right| ^{2}\leq e^{-\xi^{2}E(Z^{2})+\frac{4}{3}\left|
\xi\right| ^{3}E(\left| Z\right| ^{3})}
\end{equation*}
En particulier si $\left| \xi\right| \leq \frac{1}{2}\frac{E(Z^{2})}{E(\left|
Z\right| ^{3})}$ alors $\left| \varphi _{Z}(\xi)\right| ^{2}\leq e^{-\frac{\xi^{2}}{3}E(Z^{2})}$.
\end{lem}
\bigskip

En appliquant ce lemme on obtient ainsi

\begin{equation*}
\left| \varphi _{Y_{n,t}}\Bigl(\frac{\theta }{\sigma (t)}\Bigr)\right| ^{2}\leq 
\exp\Bigl({-\frac{\sigma _{n,t}^{2}}{\sigma ^{2}(t)}\theta ^{2}+\frac{4}{3}\frac{\left|
\theta \right| ^{3}E(\left| Y_{n,t}\right| ^{3})}{\sigma ^{3}(t)}}\Bigr)
\end{equation*}

donc

\begin{equation*}
\left| \varphi _{Z_t}(\theta )\right| ^{2}\leq \prod_{n\geq 1}
        \exp\Bigl({-\frac{\sigma _{n,t}^{2}}{\sigma ^{2}(t)}\theta ^{2}+\frac{4}{3}\frac{\left| \theta
        \right| ^{3}E(\left| Y_{n,t}\right| ^{3})}{\sigma ^{3}(t)}}\Bigr)
        = \exp\Bigl({-\theta^{2}\bigl(1-\frac{4}{3}\left| \theta \right| \frac{\Gamma _{3}(t)}{\sigma ^{3}(t)}\bigr)}\Bigr)
\end{equation*}

Pour conclure, si $\left| \theta \right| \leq 
            \displaystyle\frac{1}{4\frac{\Gamma _{3}(t)}{\sigma ^{3}(t)}}$ alors $1-\frac{4}{3}\left| \theta \right| 
            \frac{\Gamma _{3}(t)}{\sigma^{3}(t)}\geq 2/3$ 
et par cons\'equent $\left| \varphi _{Z_t}(\theta )\right| ^{2}\leq e^{-2\theta ^{2}/3}$.

$\square $

\bigskip

Comme $\pi \sigma (t)\thicksim \displaystyle \sqrt{\frac{1}{6}}\frac{\pi ^{2}}{r ^{3/2}}$ lorsque $t$ tend vers 1
il suffit maintenant d'obtenir une majoration de la fonction caract\'eristique sur le domaine 
$\displaystyle\frac{1}{4C_{3}}r ^{-1/2}\leq \left| \theta \right| \leq \sqrt{\frac{1}{6}}\frac{\pi ^{2}}{r ^{3/2}}.$

\begin{lem} \label{lemmeB}
\bf Soit $C$ une constante positive. Sous l'hypoth\`ese $\frac{C}{r ^{1/2}}\leq \left| \theta \right| <\pi\sigma (t)$ il existe un r\'eel positif $B$ tel que l'on ait $\left| \varphi _{Z_t}(\theta )\right| \leq e^{-B/r }$.
\end{lem}

\medskip

{\bf D\'emonstration.}

La m\'ethode consiste \`{a} \'ecrire

$\ln(\left| \varphi _{Z_t}(\theta )\right| )=\sum_{k\geq 1}\ln(\left|
1+t^{k}e^{ik\theta /\sigma (t)}\right| )-\ln(1+t^{k})$.

On d\'eveloppe 
\begin{equation*}
\left| 1+t^{k}e^{ik\theta /\sigma (t)}\right| ^{2}=1+t^{2k}+2t^{k}\cos
(k\theta /\sigma (t))
\end{equation*}
et on \'ecrit
\begin{eqnarray*}
\ln(\left| \varphi _{Z_t}(\theta )\right| ) &=&\frac{1}{2}\sum_{k\geq
1}\ln(1+t^{2k}+2t^{k}\cos (k\theta /\sigma (t)))-\ln(1+t^{2k}+2t^{k}) \\
&=&\frac{1}{2}\sum_{k\geq 1}\ln(1+\frac{2t^{k}(\cos (k\theta /\sigma (t))-1)}{1+t^{2k}+2t^{k}}) \\
&\leq &\frac{1}{2}\sum_{k\geq 1}\frac{2t^{k}(\cos (k\theta /\sigma (t))-1)}{1+t^{2k}+2t^{k}} \\
&\leq &\frac{1}{4}\sum_{k\geq 1}t^{k}(\cos (k\theta /\sigma (t))-1)
\end{eqnarray*}

Or on a 
\begin{equation*}
\sum_{k\geq 1}t^{k}(\cos (k\theta /\sigma (t))= {Re}\bigl(\frac{te^{i\theta
/\sigma (t)}}{1-te^{i\theta /\sigma (t)}}\bigr)=\frac{t\cos (\theta /\sigma
(t))-t^{2}}{1-2t\cos (\theta /\sigma (t))+t^{2}}
\end{equation*}

et puisque $Cr \leq \left| \theta \right| /\sigma (t)<\pi $ alors $\cos (\theta
/\sigma (t))\leq \cos (Cr )$. Par cons\'equent 
\begin{equation*}
\sum_{k\geq 1}t^{k}(\cos (k\theta /\sigma (t))\leq \frac{t\cos (Cr )-t^{2}}{1-2t\cos (Cr )+t^{2}}
\end{equation*}
donc
\begin{equation*}
\ln(\left| \varphi _{Z_t}(\theta )\right| )\leq \frac{1}{4}\Bigl(\frac{t\cos(Cr )-t^{2}}{1-2t\cos (Cr )+t^{2}}
         -\frac{t}{1-t}\Bigr)\backsim \left( \frac{1}{1+C^{2}}-1\right) r ^{-1}
\end{equation*}

On en d\'eduit l'existence d'une constante $B>0$ telle que l'on ait :
\begin{equation*}
\left| \varphi _{Z_t}(\theta )\right| \leq e^{-B/r }
\end{equation*}

$\square $

\begin{thm}
\bf On a
\begin{equation*}
\lim_{t\rightarrow 1}\int_{-\pi \sigma (t)}^{\pi \sigma (t)}\left| \varphi _{Z_t}(\theta
)-e^{-\theta ^{2}/2}\right| d\theta = 0
\end{equation*}
\end{thm}

\medskip 

{\bf D\'emonstration.}

D'apr\`es le lemme (\ref{lemmeB}) :
\begin{equation*}
\left| \varphi _{Z_t}(\theta )\right| \leq e^{-B/r }\text{ si }
                \frac{C}{r ^{1/2}}\leq \left| \theta \right| \leq \pi \sigma (t),\text{ avec }B>0
\end{equation*}

et clairement 
\begin{equation*}
e^{-\theta ^{2}/2}\leq e^{-C^{2}/2r }
\end{equation*}
sous les m\^emes conditions.
On a donc, avec $D=\min (B,C^{2}/2)$ :
\begin{eqnarray*}
\int_{\frac{C}{r ^{1/2}}\leq \left| \theta \right| \leq \pi \sigma
(t)}\left| \varphi _{Z_t}(\theta )-e^{-\theta ^{2}/2}\right| d\theta  &\leq
&\int_{\frac{C}{r ^{1/2}}\leq \left| \theta \right| \leq \pi \sigma
(t)}\left| \varphi _{Z_t}(\theta )\right| d\theta +\int_{\frac{C}{r
^{1/2}}\leq \left| \theta \right| \leq \pi \sigma (t)}e^{-\theta
^{2}/2}d\theta  \\
&\leq &e^{-D/r }\Bigl(\pi \sigma (t)-\frac{C}{r ^{1/2}}\Bigr)
\end{eqnarray*}

et cette derni\`ere quantit\'e tend vers 0 lorsque $t$ tend vers 1 ({\it i.e.\/} lorsque $r$ tend vers 0).

Il reste \`{a} voir que 
\begin{equation*}
\lim_{t\rightarrow 1}\int_{\left| \theta \right| \leq \frac{C}{r ^{1/2}}}\left| \varphi
_{Z_t}(\theta )-e^{-\theta ^{2}/2}\right| d\theta =0  
\end{equation*}

D'apr\`es le lemme (\ref{lemmeA}), sur cet intervalle on a $\left| \varphi
_{Z_t}(\theta )\right| \leq e^{-\theta ^{2}/3}$ donc 
\begin{equation*}
\left| \varphi _{Z_t}(\theta )-e^{-\theta ^{2}/2}\right| \leq e^{-\theta
^{2}/3}+e^{-\theta ^{2}/2}
\end{equation*}

et on peut conclure par le th\'eor\`eme de la convergence domin\'ee.

$\square $

\subsection{Application du th\'eor\`eme des \'equivalents}

On a choisi comme \'equivalent de la fonction $m$ lorsque $t$ tend vers 1 la fonction $m_1$ d\'efinie par

\begin{equation*}
m_1(e^{-r})=\frac{\pi^2}{12}\frac{1}{r^2}.
\end{equation*}

La d\'efinition de  $\tau_n$ par l'\'egalit\'e $m_1(\tau_n)=n$ se traduit, en posant  $\tau_n=e^{-\rho_n}$, par

\begin{equation*}
m_1(e^{-\rho_n })=n\quad\hbox{ ce qui donne }\quad \rho _{n}=\frac{1}{2\sqrt{3}\sqrt{n}}\pi
\end{equation*}

et par cons\'equent 
\begin{equation*}
\tau _{n}=e^{-\frac{1}{2\sqrt{3}\sqrt{n}}\pi}
\end{equation*}

On a aussi choisi comme \'equivalent de la fonction $\sigma$ lorsque $t$ tend vers 1 la fonction $\sigma_1$ d\'efinie par

\begin{equation*}
\sigma_1(e^{-r })=\sqrt{\frac{\pi ^{2}}{6}\frac{1}{r^3}}
\end{equation*}

ce qui donne 
\begin{equation*}
\sigma _{1}^{2}(\tau _{n})=\sigma_1^2(e^{-\rho_n })=\frac{\pi ^{2}}{6}\frac{1}{\rho_n^3}=\frac{\pi ^{2}}{6}\frac{1}{(\frac{1}{2\sqrt{3}\sqrt{n}}\pi )^{3}}
              =\frac{4}{\pi }\left( n\right) ^{3/2}\sqrt{3}.
\end{equation*}

La condition du th\'eor\`eme des \'equivalents :
\begin{equation*}
\frac{m(\tau _{n})-m_{1}(\tau _{n})}{\sigma _{1}(\tau _{n})}\rightarrow 0 \text { quand }  \tau _{n} \rightarrow 1
\end{equation*}

est bien satisfaite car on a vu en 4.1.1. que $m(e^{-r })=m_{1}(e^{-r })+O(\frac 1r),$ ce qui permet d'\'ecrire

\begin{equation*}
\frac{m(\tau _{n})-m_{1}(\tau _{n})}{\sigma _{1}(\tau _{n})}=\frac{m(e^{- \rho_{n}})-m_{1}(e^{- \rho_{n}})}{\sigma _{1}(e^{- \rho_{n}})}
=\frac{O(\frac{1}{\rho_n})}{\sigma_1(e^{-\rho_n })}=\frac{O(n^\frac 12)}
{3^{\frac 14}\sqrt{\frac{4}{\pi }} \  n^{3/4}}
\rightarrow 0
\end{equation*}

\bigskip

Le th\'eor\`eme des \'equivalents (\ref{formuleAsymptGeneralTau}) nous permet donc d'obtenir la formule asymptotique : 
$$
a_{n}\thicksim \displaystyle\frac{f(\tau _{n})}{\sqrt{2\pi }%
\sigma _{1}(\tau _{n})\tau _{n}^{n}}
$$ 

Il reste \`{a} donner un \'equivalent de
\begin{equation*}
f(\tau _{n})=\prod_{k=1}^{+\infty }(1+e^{-\frac{1}{2\sqrt{3}\sqrt{n}}\pi k})
\end{equation*}

Passons au logarithme 
\begin{equation*}
\ln (f(\tau _{n}))=\sum_{k\geq 1}\ln (1+e^{-\frac{1}{2\sqrt{3}\sqrt{n}}\pi k})
\end{equation*}

\begin{lem} 
\label{lemmeApplicEquiv}
\bf On a pour $\rho \rightarrow 0^{+}$
\begin{equation*}
\sum_{k=1}^{+\infty }\ln(1+e^{-\rho k})=\frac{\pi ^{2}}{12\rho }-\frac{1}{2}\ln 2+O\left( \rho \right)
\end{equation*}
\end{lem}
\smallskip

{\bf D\'emonstration.}

Appliquons la formule d'Euler-McLaurin: 
\begin{equation*}
\sum_{k\geq 1}\ln(1+e^{-\rho k})=\int_{1}^{+\infty }\ln (1+e^{-\rho
x})dx+\frac{1}{2}\ln (1+e^{-\rho })-\rho \int_{1}^{+\infty }b_{1}(x)\frac{e^{-\rho x}}{1+e^{-\rho x}}dx
\end{equation*}

{\bf a) Le terme $\int_{1}^{+\infty }\ln (1+e^{-\rho x})dx$ :}
\par\smallskip
On d\'ecompose l'int\'egrale:

\begin{eqnarray*}
\int_{1}^{+\infty }\ln (1+e^{-\rho x})dx &=&\int_{0}^{\infty }\ln
(1+e^{-\rho x})dx-\int_{0}^{1}\ln (1+e^{-\rho x})dx \\
&=&\int_{0}^{+\infty }\sum_{n\geq 1}\frac{(-1)^{n-1}}{n}e^{-\rho
nx}dx-\ln 2+O\left( \rho \right) \\
&=&\sum_{n\geq 1}\frac{(-1)^{n-1}}{\rho n^{2}}-\ln 2+O\left(\rho\right)
\end{eqnarray*}

Donc 
\begin{equation*}
\int_{1}^{+\infty }\ln (1+e^{-\rho x})dx=-\ln 2+\frac{\pi ^{2}}{12\rho }+O\left(\rho \right)
\end{equation*}

{\bf b) Le terme $\frac{1}{2}\ln (1+e^{-\rho})$ :}
\par\smallskip
\begin{equation*}
\frac{1}{2}\ln (1+e^{-\rho })=\frac{1}{2}\ln 2+O\left( \rho \right)
\end{equation*}

{\bf c) Le troisi\`eme terme :}
\par\smallskip
On a 
\begin{equation*}
-\rho \int_{1}^{+\infty }b_{1}(x)\frac{e^{-\rho x}}{1+e^{-\rho x}}dx=O\left(\rho e^{-\rho }\right)
\end{equation*}

En effet la fonction $x\mapsto \displaystyle\frac{e^{-\rho x}}{1+e^{-\rho x}}$ est
positive d\'ecroissante et elle tend vers 0 \`{a} l'infini. Comme la fonction $b_{1}$
est p\'eriodique, par le lemme d'Abel on obtient la majoration 
\begin{equation*}
\left| \int_{1}^{+\infty }b_{1}(x)\frac{e^{-\rho x}}{1+e^{-\rho x}}dx\right| 
               \leq C\frac{e^{-\rho }}{1+e^{-\rho }}\leq Ce^{-\rho }
\end{equation*}

$\square $

\textbf{Conclusion}

D'apr\`es le lemme  (\ref{lemmeApplicEquiv}) avec $\rho=\rho_n$ on a 
\begin{equation*}
\sum_{k\geq 1}\ln (1+e^{-\frac{1}{2\sqrt{3}\sqrt{n}}\pi k})
           = \frac{\pi \sqrt{n}}{2\sqrt{3}}-\frac{1}{2}\ln 2+O\left( \frac{1}{\sqrt{n}}\right)
\end{equation*}

donc

\begin{equation*}
f(\tau _{n})=\prod_{k=1}^{+\infty }(1+e^{-\frac{1}{2\sqrt{3}\sqrt{n}}\pi
k})\thicksim \frac{1}{\sqrt{2}}e^{\frac{\pi \sqrt{n}}{2\sqrt{3}}}
\end{equation*}
lorsque $n$ tend vers l'infini, ce qui donne la formule asymptotique des partitions restreintes :
\begin{equation*}
q(n)\thicksim \frac{1}{\sqrt{2}}e^{\frac{\pi \sqrt{n}}{2\sqrt{3}}}
               \frac{1}{\sqrt{2\pi }\sqrt{\frac{4}{\pi }\left( n\right) ^{3/2}\sqrt{3}}
               e^{-\frac{\sqrt{n}}{2\sqrt{3}}\pi }}=\frac{1}{4}\frac{e^{\frac{\pi \sqrt{n}}{\sqrt{3}}}}{3^{1/4}n^{3/4}}
\end{equation*}

\bigskip

\section{Appendice}

\subsection{D\'emonstration du Lemme 3.3.}

\smallskip
Comme $\sup_{n\geq 1}|u_{n,t}|\rightarrow 0$ quand $t\rightarrow 1,$ il
existe $a<1$ tel que  pour $t\in\ ]a,1[$ on a  $|u_{n,t}|<1/2$ 
pour tout $n\geq 1.$ Donc $\ln(1+u_{n,t})$ est bien d\'{e}fini pour $t\in\ ]a,1[$ et 
\begin{equation*}
\ln(1+u_{n,t})=\sum_{k=1}^{+\infty }\frac{(-1)^{k-1}}{k}(u_{n,t})^{k}
\end{equation*}

ce qui donne  
\begin{equation*}
\left| \ln(1+u_{n,t})-u_{n,t}\right| \leq \left| u_{n,t}\right|
^{2}\sum_{k=2}^{+\infty }\frac{1}{k}\left| u_{n,t}\right| ^{k-2}\leq \left|
u_{n,t}\right| ^{2}\sum_{k=2}^{+\infty }(\frac{1}{2})^{k-2}=2\left|
u_{n,t}\right| ^{2}
\end{equation*}

D'autre part la s\'{e}rie  $\sum_{n\geq 1}|u_{n,t}|$ est suppos\'{e}e
convergente pour tout $t\in\ ]\alpha ,1[,$ donc la s\'{e}rie $\sum_{n\geq
1}\left| u_{n,t}\right| ^{2}$ est convergente si $t\in\ ]\sup (a,\alpha ),1[.$
On en d\'{e}duit que la s\'{e}rie $\sum_{n\geq 1}\ln(1+u_{n,t})$ est
convergente si $t\in\ ]\sup (a,\alpha ),1[$ et il en est donc de m\^{e}me du
produit infini $\prod_{n\geq 1}(1+u_{n,t}).$ 

\bigskip 

D'autre part, pour tout $N\geq 1$ on a  
\begin{equation*}
\left| \sum_{n=1}^{N}\ln(1+u_{n,t})-\sum_{n=1}^{N}u_{n,t}\right| \leq
\sum_{n=1}^{N}\left| \ln(1+u_{n,t})-u_{n,t}\right| \leq
2\sum_{n=1}^{N}\left| u_{n,t}\right| ^{2}
\end{equation*}

Comme 
\begin{equation*}
\sum_{n=1}^{N}\left| u_{n,t}\right| ^{2}\leq \sup_{n\geq
1}|u_{n,t}|\sum_{n=1}^{N}\left| u_{n,t}\right| \leq M\sup_{n\geq
1}|u_{n,t}|
\end{equation*}

on en d\'{e}duit que $\sum_{n=1}^{+\infty }\left| u_{n,t}\right| ^{2}\leq M\sup_{n\geq 1}|u_{n,t}|$ et  que  
\begin{equation*}
\lim_{N\rightarrow +\infty }\left|
\sum_{n=1}^{N}\ln(1+u_{n,t})-\sum_{n=1}^{N}u_{n,t}\right| \leq 2M\sup_{n\geq 1}|u_{n,t}|
\end{equation*}

Donc  
\begin{equation*}
\left| \sum_{n=1}^{+\infty }\ln(1+u_{n,t})-\sum_{n=1}^{+\infty
}u_{n,t}\right| \leq 2M\sup_{n\geq 1}|u_{n,t}|.
\end{equation*}

Pour conclure il suffit de prendre la limite quand $t\rightarrow 1$, on
obtient 
\begin{equation*}
\lim_{t\rightarrow 1}\sum_{n=1}^{+\infty }\ln(1+u_{n,t})=\lim_{t\rightarrow
1}\sum_{n=1}^{+\infty }u_{n,t}=S
\end{equation*}

En passant \`{a} l'exponentielle on obtient 
\begin{equation*}
\lim_{t\rightarrow 1}\prod_{n=1}^{+\infty }(1+u_{n,t})=e^{S}
\end{equation*}

$\square $

\subsection{D\'emonstration du Lemme de Cram\'er (Lemme \ref{lemmeCramer})} \label{preuveLemCramer} 
{\bf D\'emonstration.} ({\it voir\/} Cram\'er~\cite{Cramer}, Chung~\cite{Chung} p. 210)
\par\smallskip
Soit $Y$ une variable al\'eatoire ind\'{e}pendante de $Z$ et de m\^{e}me loi. On a
\begin{equation*}
\left| \varphi _{Z}(\xi)\right| ^{2}=\varphi _{Z}(\xi)\overline{\varphi _{Y}(\xi)}
          = E(e^{i\xi Z})\overline{E(e^{i\xi Y})}=E(e^{i\xi (Z-Y)})
\end{equation*}

ce qui permet d'\'{e}crire
\begin{equation*}
\left| \varphi _{Z}(\xi)\right| ^{2}=\int_{\mathbb{R}^{2}}e^{i\xi (z-y)}dP_{Z}(z)dP_{Y}(y)
               =\int_{\mathbb{R}^{2}}\cos(\xi (z-y))dP_{Z}(z)dP_{Y}(y)
\end{equation*}

En utilisant la majoration 
\begin{equation*}
\cos (u)\leq 1-\frac{u^{2}}{2}+\frac{\left| u\right| ^{3}}{6}
\end{equation*}
que l'on peut obtenir \`a l'aide de la formule de Taylor d'ordre deux avec reste int\'egral,
on en d\'{e}duit que 
\begin{equation*}
\left| \varphi _{Z}(\xi)\right| ^{2}\leq 1-\frac{\xi^{2}}{2}\int_{\mathbb{R}^{2}}(z-y)^{2}dP_{Z}(z)dP_{Y}(y)
           +\frac{\left| \xi\right| ^{3}}{6}\int_{\mathbb{R}^{2}}\left| z-y\right| ^{3}dP_{Z}(z)dP_{Y}(y)
\end{equation*}

La premi\`{e}re int\'{e}grale n'est autre que $2E(Z^{2}).$ Pour  la deuxi\`{e}me  on utilise la majoration
\begin{equation*}
\left| z-y\right| ^{3}\leq 4\left| z\right| ^{3}+4\left| y\right| ^{3}
\end{equation*}

ce qui permet de majorer l'int\'{e}grale par $8E(\left| Z\right| ^{3}).$

On obtient finalement
\begin{equation*}
\left| \varphi _{Z}(\xi)\right| ^{2}\leq 1-\xi^{2}E(Z^{2})+\frac{4}{3}\left|
\xi\right| ^{3}E(\left| Z\right| ^{3})\leq e^{-\xi^{2}E(Z^{2})+\frac{4}{3}\left|
\xi\right| ^{3}E(\left| Z\right| ^{3})}
\end{equation*}

Pour la seconde partie du lemme, si $\left| \xi\right| \leq \frac{1}{2}\frac{E(X^{2})}{E(\left| X\right| ^{3})}$ 
il suffit de remarquer que  
\begin{equation*}
-\xi^{2}E(Z^{2})+\frac{4}{3}\left| \xi\right| ^{3}E(\left| Z\right|
^{3})=-\xi^{2}[E(Z^{2})-\frac{4}{3}\left| \xi\right| E(\left| Z\right|
^{3})]\leq -\xi^{2}\frac{1}{3}E(Z^{2})
\end{equation*}

\par\smallskip
$\square$

\medskip

\textbf{Remerciements}

Nous remercions Mesdames C\'ecile Fouilh\'e et No\" emie El Qotbi pour l'int\'er\^et qu'elles ont port\'e \`a l'\'etude de l'article de Luis B\' aez-Duarte.


\begin{thebibliography}{9}

\bibitem{LBD} B\' aez-Duarte, L., 
   {\it ``Hardy-Ramanujan's Asymptotic Formula for Partitions and the Central Limit Theorem,''}
          Advances in Mathematics 125, 114-120 (1997). 
\bibitem{Chung} Chung, K.L., 
   {\it ``A Course in Probability Theory, 3rd ed.''}
          Academic Press, San Diego [CA] (2001). 
\bibitem{Cramer} Cram\'er, H., 
   {\it ``Random Variables and Probability Distributions,''}
          2nd ed., Cambridge Univ. Press, Cambridge (1963). 
\bibitem{Erdos} Erdös, P., Lehner, J.,
   {\it ``The Distribution of the Number of Summands in the Partitions of a Positive Integer,''}
      Duke Mathematical Journal Vol. 8, No.2 (June, 1941)          
\bibitem{Hayman} Hayman, W.K., 
   {\it ``A Generalisation of Stirling's Formula,''}
             J. Reine Angew. Mat. 196, Nos. 1/2, 67-95 (1956).
\bibitem{Ingham} Ingham, A.E., 
   {\it ``A Tauberian Theorem for Partitions,''}
          The Annals of Mathematics, Second Series, Vol. 42, No. 5, 1075-1090 (Dec., 1941)
\bibitem{Rosen} Rosenbloom, P.C., 
   {\it ``Probability and Entire Functions,''}
          Studies in Mathematical Analysis and Related Topics, 
          Vol. 45, 325-332, Stanford Univ. Press, Palo Alto, CA (1962).      
          
\end{thebibliography}
\end{document}